%
%
%

\input amstex.tex
\input amsppt.sty
\input pictex.tex
\magnification=\magstep1
\advance\vsize-0.5cm\voffset=-0.5cm\advance\hsize1cm\hoffset0cm
\NoBlackBoxes

\def\delet{\mathaccent"7017 }

\def\Int{\mathop{\fam0 Int}}

\def\R{\Bbb R} \def\Z{\Bbb Z} 

\topmatter
\title Codimension two PL embeddings of spheres with
nonstandard regular neighborhoods
\endtitle
\rightheadtext{Codimension two PL embeddings}
\author M. Cencelj, D. Repov\v s and A. Skopenkov \endauthor
\address Institute of Mathematics, Physics and Mechanics,
University of Ljubljana, P.~O.~Box 2964, 1001 Ljubljana, Slovenia.
e-mail: dusan.repovs\@fmf.uni-lj.si \endaddress
\address Institute of Mathematics, Physics and Mechanics,
University of Ljubljana, P.~O.~Box 2964, 1001 Ljubljana, Slovenia.
e-mail: matija.cencelj\@fmf.uni-lj.si \endaddress
\address Department of Differential Geometry, Faculty of Mechanics and
Mathematics, Moscow State University, Moscow, Russia 119992.
e-mail: skopenko\@mccme.ru \endaddress
\thanks
Skopenkov gratefully acknowledges the support of the Pierre Deligne fund based on
his 2004 Balzan prize in mathematics, INTAS Grant No. YSF-2002-393, the Russian Foundation for
Basic Research, Grants No 05-01-00993 and 04-01-00682, President of Russian Federation Grants
MD-3938.2005.1 and NSH-1988.2003.1. Cencelj and Repov\v{s} were supported by Slovenian Research Agency grants BI-RU/05-07-04 and BI-RU/05-07-13.
\endthanks
\abstract 
For a given polyhedron $K\subset M$ the notation $R_M(K)$ denotes a regular
neighborhood of $K$ in $M$.
We study the following problem: find all pairs $(m,k)$ such that if $K$ is a
compact $k$-polyhedron and $M$ a PL $m$-manifold, then
$R_M(fK)\cong R_M(gK)$, for each two homotopic PL embeddings $f,g:K\to M$.
We prove that $R_{S^{k+2}}(S^k)\not\cong S^k\times D^2$ for each $k\ge2$ and
{\it some} PL sphere $S^k\subset S^{k+2}$  (even for {\it any} PL sphere
$S^k\subset S^{k+2}$ having an isolated non-locally flat point with the
singularity $S^{k-1}\subset S^{k+1}$ such that
$\pi_1(S^{k+1}-S^{k-1})\not\cong\Z$).
\endabstract 
\endtopmatter

\document

This paper is on the uniqueness of {\it regular neighborhoods} problem 
for distinct embeddings of one manifold into the other. 
Definitions of regular neighborhood and other notions of PL topology can 
be found in [RoSa72]. 
For any subpolyhedron $X\subset M$ of a PL manifold $M$ we denote  the
regular neighborhood of $X$ in $M$ by $R_M(X)$.
The following is the main result of the present paper:

\proclaim{Theorem 1} $R_{S^{k+2}}(S^k)$ is not homeomorphic to 
$S^k\times D^2$ for any $k\ge2$ and

(a) any PL sphere $S^k\subset S^{k+2}$ which is the suspension over a locally
flat PL knot $S^{k-1}\subset S^{k+1}$ such that
$G = \pi_1(S^{k+1}-S^{k-1})\not\cong\Z$; or

(b) any PL sphere $S^k\subset S^{k+2}$ having an isolated non-locally flat
point with the singularity $S^{k-1}\subset S^{k+1}$ such that
$G = \pi_1(S^{k+1}-S^{k-1})\not\cong\Z$.
\endproclaim

For $k=2$ in Theorem 1.b one can take {\it any} non-locally flat PL sphere 
$S^2\subset S^4$. 

Recall that for $x\in S^k\subset S^{k+2}$ {\it the singularity at point $x$}
is the isotopy class of the submanifold $\partial D^{k+2}\cap S^k$ of
$\partial D^{k+2}$. Here $(D^{k+2},D^{k+2}\cap S^k)$ is any PL $(k+2,k)$ ball 
pair which is a regular neighborhood of $x$ in the pair $(S^{k+2},S^k)$. 

Of course, Theorem 1.a follows from Theorem 1.b, but the proof of the
former is simpler (in particular for $k=2$, because our argument uses the
fundamental group whereas the proof of Theorem 1.b for that cases uses [Gab86, Gab87]). 

We now give some motivation for Theorem 1.1.
An important question in the topology of manifolds [Mas59] asks for which cases
the normal vector bundle of a smooth $k$-manifold $K$ is independent of the
choice of the smooth embedding $K\to S^m$.  This is the case when:

(a) $m=k+1$; or

(b) $m=k+2$ and $K$ is orientable [Ker65], [Rol76]; or

(c) $K$ is a homotopy $k$-sphere and either $m=k+3$ or $m\ge\frac{3k}2+1$
[Ker59], [Mas59], [HaHi62]; or

(d) $K=\R P^2$ and $m=4$ [Mas69].

On the other hand, there exist smooth embeddings $K\subset S^m$ with
{\it distinct} normal bundles [HLS65; Haefliger's example in \S2], [Mas69],
[Roh70], [MiRe71].
For example,  $K$ can be the Klein bottle and $f_0,f_4:K\to S^4$ can be any
smooth embedings with the normal Euler classes equal to 0 and 4, respectively.
Then $\pi_1(R_{S^4}(f_0K))\cong\Z^3$ while $\pi_1(R_{S^4}(f_4K))$ is not
virtually abelian, i.e. it does not contain any abelian subgroups of finite
index.
This argument, due to Hillman, shows that for such embeddings even the
{\it spaces} of normal bundles are {\it non-homeomorphic}.

Any locally flat PL embedding has a normal block bundle [RoSa68].
Thus the Massey question makes sense in the locally flat PL category.
In this category the normal bundle does not depend on the embedding when
$K=S^k$ (by the Zeeman Unknotting Balls Theorem) and in the cases (a), (b),
(d) above [Wal67; Theorem 2], [RoSa68; Corollary 6.8].
In particular, $R_{S^4}(S^2)\cong S^2\times D^2$ for any locally flat PL
embedding $S^2\subset S^4$ (cf. Theorem 1).
Note that a PL sphere $S^n\subset\R^m$ has a PL standard regular neighborhood 
if either $m-n\ne2$ or the embedding $S^n\subset\R^m$ is locally flat [Wal67, 
Zee63]. 

From now on we shall work in the PL (not locally flat PL) category, unless
otherwise specified.
As another polyhedral version of the Massey problem we can consider the
question, whether for embeddings $f,g:K\to S^m$ the homeomorphism
$g\circ f^{-1}:fK\to gK$ can be extended to a homeomorphism
$R_M(fK)\cong R_M(gK)$ [LiSi69].
The analogues of Theorem 1 and conjectures (e), (f), (g) below for this
problem are trivial.
Theorem 1 was motivated by the following related problem, raised
independently by Shtanko and Chernavski (private communication) and also in 
[Hor85, Hor90, Rep88]:

{\it Find all pairs $(m,k)$ such that if $K$ is a compact $k$-polyhedron 
and $M$ a PL $m$-manifold, then

(*) $R_M(fK)\cong R_M(gK)$ for each two homotopic PL embeddings $f,g:K\to M$.}

The particular case of this problem for $M=S^m$ (then $f$ and $g$ are always
homotopic) attracted special attention.

The property (*) holds for $m\ge2k+2$ by general position.
In general, many invariants of $R_M(fK)$ and $R_M(gK)$ coincide: the homotopy
types, the intersection rings, the higher Massey products, the (classifying
maps of) tangent bundles (and hence also all invariants deduced from the
tangent bundle, e.g. characteristic classes and numbers).
This implies that (*) also holds for $m=2k+1$ (because in this case an
$m$-thickening is completely determined by its tangent bundle [LiSi69]) and
for $m=2$ (because a 2-manifold $N$ with boundary is completely determined by
$H_1(N,\Z_2)$ and its intersection form).
Also,
$$R_M(fK)\times I^l\cong R_{M\times I^l}(fK)\cong R_{M\times I^l}(gK)
\cong R_M(gK)\times I^l$$
for $l\ge2k+1-m$.
The boundary $\partial R_M(fK)$ is $l$-connected for $m\ge k+l+2$ and
$l$-connected $K$.
All this suggests that distinguishing $R_M(fK)$ and $R_M(gK)$ is a
nontrivial problem for $m\ge k+3$.

The property (*) holds for $m=k+1\ge3$ and a {\it fake surface} $K$
[Cas65], [Cav92], [Mat73], [Rep88], [RBS99].
The property (*) fails for:

(a) $M=S^m$, $m=k+1\ge3$ and $K=S^k\vee S^1\vee S^1$ [Rep88]; or

(b) $M=S^3$, $K$ a common spine of the granny knot and the square knot
[MPR89], [CaHe94], [CLR97], [CHR98]; or

(c) $M=S^4$ and the Dunce Hat $K$ [Zee64], [Hor85], [Hor90]; or

(d) $M=S^{k+2}$, $k\ge2$ and $K=S^k$ by Theorem 1.

We conjecture that the property (*) also fails for:

(e) $M=S^{2k}$ and $K=\Sigma(S^{k-1}\sqcup S^{k-1})$ 
(one can take $f$ and $g$ to be the suspensions over the trivial link 
and the Hopf link, respectively); or

(f) $M=S^4$ and $K=\R P^2$ ($f$ and $g$ should be non-locally flat); or

(g) each $k+3\le m\le2k$ and some polyhedron or even manifold $K$.

The case (a) above was proved using the number of connected components of
$\partial R_M(fK)$; the cases (c) and (d) was and will be proved using
$\pi_1(\partial R_M(fK))$.

\smallskip
{\it Proof of Theorem 1.}
If a PL embedding $S^2\subset S^4$ is locally flat everywhere except at $n$ 
points with singularities $L_1,\dots,L_n$, then $\partial R_{S^4}(S^2)$ is 
obtained from $\partial D^4$ by the Dehn surgery over the knot  
$L_1\#\dots\#L_n$  with certain framing (the details are analogous to the proof 
of Theorem 2.a below). 
Hence the case $k=2$ of Theorem 1.b follows by a deep result of Gabai 
[Gab86, Gab87; Remark after Corollary 5], to the effect that 
{\it the Dehn surgery along any non-trivial knot $S^1\subset S^3$ can never 
yield $S^1\times S^2$}.
(The argument of [Gab86, Gab87] uses foliation theory and cannot be generalized 
to obtain Theorem 1.b for $k>2$.)

Theorem 1.b for $k>2$ and Theorem 1.a are implied by the following result. 
\qed

\proclaim{Theorem 2} For $k\ge2$ and given embedding $S^k\subset S^{k+2}$ 
set $N=R_{S^{k+2}}(S^k)$.
Then $G:=\pi_1(S^{k+1}-S^{k-1})$ maps into $\pi_1(\partial N)$ in each of 
the following two cases: 

(a) the embedding $S^k\subset S^{k+2}$ is the suspension over a locally
flat PL knot $S^{k-1}\subset S^{k+1}$. 

(b) $k\ge3$ and the embedding $S^k\subset S^{k+2}$ has an isolated non-locally 
flat point with the singularity $S^{k-1}\subset S^{k+1}$.
\endproclaim

\beginpicture 
\setcoordinatesystem units <1.000mm,1.000mm>
\setplotarea x from -40 to 40, y from -50 to 40
\put {Figure 1.a} at 0 -45
\setsolid
\setplotsymbol ({\fiverm.})
\plot -30.000 0.000 30.000 0.000 /
\arrow <2mm> [.2,.4] from -30.000 10.000 to -15.000 0.000
\arrow <2mm> [.2,.4] from -30.000 10.000 to 15.000 0.000
\plot -10.000 0.000 -4.000 -12.000 /
\plot 4.000 -12.000 10.000 0.000 /
\plot -20.000 0.000 -4.000 -32.000 /
\plot 4.000 -32.000 20.000 0.000 /
\circulararc 126.870 degrees from -4.000 -12.000 center at 0.000 -10.000
\circulararc 126.870 degrees from -4.000 -32.000 center at -0.000 -30.000
\arrow <2mm> [.2,.4] from 10.000 -30.000 to 5.000 -25.000
\arrow <2mm> [.2,.4] from 15.000 -15.000 to 10.000 -10.000
\arrow <2mm> [.2,.4] from -20.000 -30.000 to -5.000 -10.000
\arrow <2mm> [.2,.4] from -20.000 -30.000 to -15.000 -10.000
\arrow <2mm> [.2,.4] from -20.000 30.000 to -15.000 10.000
\arrow <2mm> [.2,.4] from -20.000 30.000 to -5.000 10.000
\arrow <2mm> [.2,.4] from 15.000 15.000 to 10.000 10.000
\arrow <2mm> [.2,.4] from 10.000 30.000 to 5.000 25.000
\circulararc 126.870 degrees from 4.000 32.000 center at 0.000 30.000
\circulararc 126.870 degrees from 4.000 12.000 center at -0.000 10.000
\plot 4.000 32.000 20.000 0.000 /
\plot -20.000 0.000 -4.000 32.000 /
\plot 4.000 12.000 10.000 0.000 /
\plot -10.000 0.000 -4.000 12.000 /
\setplotsymbol ({\rm.})
\plot -15.000 0.000 0.000 -30.000 /
\plot 0.000 -30.000 15.000 0.000 /
\plot 0.000 30.000 15.000 0.000 /
\plot -15.000 0.000 0.000 30.000 /
\put {$D_-^{k+2}$} at 40.000 -30.000
\put {$D_+^{k+2}$} at 40.000 30.000
\put {$S^{k-1}$} at -35.000 10.000
\put {$N_-$} at 15.000 -35.000
\put {$B_-^k$} at 20.000 -20.000
\put {$\partial_-N$} at -20.000 -35.000
\put {$\partial_+N$} at -20.000 35.000
\put {$B_+^k$} at 20.000 20.000
\put {$N_+$} at 15.000 35.000
\put {$\bullet$} at -15.000 0.000
\put {$\bullet$} at 15.000 -0.000

\setshadesymbol ({\fiverm.})
\setshadegrid span <0.02in>
\vshade -20 0 0 <z,z,,> -10  0 20
                <z,z,,>  -4 12 32
                <z,z,,>  -2 14 34
                <z,z,,>   0 14.5 34.5
                <z,z,,>   2 14 34
                <z,z,,>   4 12 32
                <z,z,,>  10  0 20
                <z,z,,>  20  0  0 /
\vshade -20 0 0 <z,z,,> -10 -20 0
                <z,z,,>  -4 -32 -12
                <z,z,,>  -2 -34 -14
                <z,z,,>   0 -34.5 -14.5
                <z,z,,>   2 -34 -14
                <z,z,,>   4 -32 -12
                <z,z,,>  10 -20 0
                <z,z,,>  20  0  0 /

\endpicture

{\it Proof of Theorem 2.a.} 
Consider decomposition
$S^{k+2}=D^{k+2}_+\bigcup\limits_{\partial D^{k+2}_\pm=S^{k+1}} D^{k+2}_-$.
Let $B^k_\pm:=D^{k+2}_\pm\cap S^k$ be the cone over $S^{k-1}$ (Figure 1.a).
Denote
$$N_\pm:=R_{D^{k+2}_\pm}(B^k_\pm)\quad\text{so that}
\quad N_\pm\cap\partial D^{k+2}_\pm=R_{S^{k+1}}(S^{k-1})\quad\text{and}$$
$$\partial_\pm N:=\partial N_\pm-\Int R_{\partial N_\pm}(S^{k-1})\quad
\text{so that}\quad N=N_+\cup N_-\quad\text{and}
\quad\partial N=\partial_+N\cup\partial_-N.$$
Since $B^k_\pm$ is a cone, it follows that $D^{k+2}_\pm$ collapses to 
$B^k_\pm$. 
Therefore $D^{k+2}_\pm$ is also a regular neighborhood of $B^k_\pm$ in
$D^{k+2}_\pm$.
Hence by the uniqueness of regular neighborhoods,
$$(\partial N_\pm,S^{k-1})\cong(S^{k+1},S^{k-1}),\quad\text{thus}\quad
\pi_1(\partial_\pm N)\cong G.$$
The composition
$\partial_+N\overset{i}\to\to\partial_+N\cup\partial_-N\overset{r}\to\to
\partial_+N$ of
the inclusion $i$ and the 'symmetric' retraction $r$ is a homeomorphism.
So the induced composition 
$G\overset{i_*}\to\to\pi_1(\partial N)\overset{r_*}\to\to G$ is an
isomorphism. Hence $i_*$ is a monomorphism, and (a) is proved.
\qed

\beginpicture 
\setcoordinatesystem units <1.000mm,1.000mm>
\setplotarea x from -50 to 50, y from -60 to 50


\setsolid
\setplotsymbol ({\fiverm.})
\arrow <2mm> [.2,.4] from -30.000 -5.000 to 0.000 0.000
\arrow <2mm> [.2,.4] from -30.000 -5.000 to -25.000 0.000
\arrow <2mm> [.2,.4] from -30.000 -30.000 to -10.000 -20.000
\arrow <2mm> [.2,.4] from -30.000 -30.000 to -20.000 -30.000
\arrow <2mm> [.2,.4] from 20.000 -45.000 to 10.000 -35.000
\arrow <2mm> [.2,.4] from 25.000 -30.000 to 15.000 -20.000
\arrow <2mm> [.2,.4] from 15.000 15.000 to 10.000 10.000
\plot -30.000 0.000 30.000 0.000 /
\plot 20.000 0.000 20.000 -30.000 /
\plot -20.000 0.000 -20.000 -30.000 /
\circulararc 120.000 degrees from -20.000 -30.000 center at -10.000 -30.000
\circulararc 120.000 degrees from 5.000 -38.660 center at 10.000 -30.000
\circulararc 60.000 degrees from 5.000 -38.660 center at 0.000 -47.321
\plot 10.000 0.000 10.000 -20.000 /
\plot -10.000 -20.000 -10.000 0.000 /
\circulararc 120.000 degrees from -10.000 -20.000 center at -5.000 -20.000
\circulararc 120.000 degrees from 2.500 -24.330 center at 5.000 -20.000
\circulararc 60.000 degrees from 2.500 -24.330 center at 0.000 -28.660
\setplotsymbol ({\rm.})
\circulararc 120.000 degrees from -15.000 -25.000 center at -7.500 -25.000
\circulararc 120.000 degrees from 3.750 -31.495 center at 7.500 -25.000
\circulararc 60.000 degrees from 3.750 -31.495 center at 0.000 -37.990
\plot -15.000 -25.000 -15.000 0.000 /
\plot -15.000 0.000 0.000 30.000 /
\plot 0.000 30.000 15.000 0.000 /
\plot 15.000 0.000 15.000 -25.000 /
\put {$N_+=D^{k+2}_+$} at 40.000 25.000
\put {$\partial_+N$} at -35.000 -5.000
\put {$\partial_-N$} at -35.000 -30.000
\put {$N_-$} at 25.000 -50.000
\put {$B_-^k$} at 30.000 -35.000
\put {$B_+^k$} at 15.000 20.000
\put {$\bullet$} at 15.000 0.000
\put {$\bullet$} at -15.000 0.000
\put {Figure 1.b} at 0 -55

\setshadesymbol ({\fiverm.})
\setshadegrid span <0.02in>
\vshade -30 0 40 <z,z,,> 30 0 40 /
\vshade -20 -30 0 <z,z,,> -18.4 -36.3 0
                  <z,z,,> -15 -38.8 0
                  <z,z,,> -10 -40 0 /
\vshade -10 -40 -20 <z,z,,>  -9 -40 -23.7
                    <z,z,,>  -5 -39 -25.2
                    <z,z,,>  -3 -37.7 -24.5
                    <z,z,,>   0 -37.2 -23.7
                    <z,z,,>   3 -37.7 -24.5
                    <z,z,,>   5 -39 -25.2
                    <z,z,,>   9 -40 -23.7
                    <z,z,,>  10 -40 -20 /
\vshade 10 -40 0 <z,z,,> 15 -38.8 0
                 <z,z,,> 18.4 -36.3 0
                 <z,z,,> 20 -30 0 /

\endpicture

{\it Proof of Theorem 2.b.} 
Let $x$ be the isolated non-locally flat point.
Then the singularity $S^{k-1}\subset S^{k+1}$ is locally flat.
Consider the decomposition of $S^{k+2}$ as above such that
$x\in\delet D^{k+2}_+$.
Let $B^k_+:=D^{k+2}_+\cap S^k$ be the cone over $S^{k-1}$ with the vertex $x$
(Figure 1.b).
Let $B^k_-:=D^{k+2}_-\cap S^k$ and $N_+:=D^{k+2}_+$.
Define $N_-$ and $\partial_\pm N$ as in the proof of (a).
Then 
$$N\cong N_+\cup N_-\quad\text{and}\quad\partial N=\partial_+N\cup\partial_-N.$$ 
We have
$$\pi_1(\partial_+N)\cong G\quad\text{and}\quad
\pi_1(\partial_+N\cap\partial_-N)\cong\pi_1(S^1\times\partial D^k)\cong\Z
\quad\text{for}\quad k\ge3.$$

First consider the simple case when $x$ is the {\it only} non-locally flat 
point (formally, the proof of the general case does not use the simpler case). 
Then $B^k_-$ is locally flat.
Since $N_-=R_{D^{k+2}_-}(B^k_-)$, it follows that the pair
$(N_-,B^k_-)$ is standard [Zee63].
Therefore $\partial_-N\cong S^1\times D^k$.
Hence $\pi_1(\partial N)\cong G$ by the van Kampen Theorem.  

In the general case denote by $a$ the generator of
$\pi_1(\partial_+N\cap\partial_-N)\cong H_1(\partial_+N\cap\partial_-N;\Z)\cong
\Z$. 
Observe that $a$ is represented by a small circle in
$\partial R_{\partial N_\pm}(S^{k-1})$, bounding a small 2-disk in
$R_{\partial N_\pm}(S^{k-1})$, transversal to $S^{k-1}$.
By [Hud69; p. 57], $\partial N_-\cong S^{k+1}$.
Hence the class of $a$ is the generator of $H_1(\partial_\pm N,\Z)\cong\Z$.
Thus the maps of $H_1(\cdot;\Z)$ induced by the inclusions 
$\partial_+N\cap\partial_-N\to\partial_\pm N$  are injective. 
Therefore the maps of $\pi_1$ unduced by the inclusions  are also injective.
Hence by the van Kampen Theorem [Rol76, pages 370-371] the inclusion-induced 
map $\pi_1(\partial_+N)\to\pi_1(\partial N)$ is also injective.
\qed

\smallskip
We remark that under the assumptions of Theorem 2 the group $\pi_1(\partial N)$ 
maps onto $G:=\pi_1(S^{k+1}-S^{k-1})$. 

We remark that for $k\ge3$ there is a PL embedding $S^k\subset S^{k+2}$ which 
is locally flat {\it except at one point} and such that 
$R_{S^{k+2}}(S^k)\not\cong S^k\times D^2$. 
Indeed, there exist slice knots $S^{k-1}\subset S^{k+1}$ having non-cyclic
$\pi_1(S^{k+1}-S^{k-1})$ [Ker65], and for each slice knot
$S^{k-1}\subset S^{k+1}$ there is an embedding $S^k\subset S^{k+2}$, having
only one non locally-flat point with the singularity $S^{k-1}\subset S^{k+1}$
[FoMi66].
Now the remark follows by (the simple case of) Theorem 2.b. 
The remark also follows from Theorem 2.a 
because the regular neighborhood of the suspension over a knot 
$L$ homeomorphic to that the regular neighborhood of certain knot with the 
only singular point (having singularity $L\#(-L)$). 
Or else cf. [Gab86, Gab87; Remark after Corollary 5].

Theorem 1 was announced in [ORS01] and at the International Conference on
Knots, Links and Manifold (Siegen, January 2001).
We would like to acknowledge J. Hillman, A. A. Klyachko, W. B. R. Lickorish,
E. Rafikov, M. A. Shtanko, J. Vrabec, and the referee for 
comments and suggestions.

\enddocument
\end